\documentclass[11pt]{amsart}
\usepackage{latexsym}
\usepackage{amscd}
\usepackage{amssymb}
\usepackage{amsmath}
\usepackage{amsthm}
\usepackage{enumerate}
\usepackage{verbatim}
\def\co{\colon\thinspace}
\def\d{\delta}
\def\e{\epsilon}

\def\g{\gamma}
\newtheorem{thm}{Theorem}[section]
\newtheorem{cor}[thm]{Corollary}

\newtheorem{lem}[thm]{Lemma}
\newtheorem{prop}[thm]{Proposition}

\newtheorem{Example}[thm]{Example}
\newenvironment{ex}{\begin{Example}\rm}{\end{Example}}
\newtheorem{remark}[thm]{Remark}
\newenvironment{rmk}{\begin{remark}\rm}{\end{remark}}
\newtheorem{Fact}[thm]{Fact}

\newtheorem{Main Observation}[thm]{Main Observation}
\newenvironment{main observation}
{\begin{Main Observation}\rm}{\end{Main Observation}}
\newtheorem{Convention}[thm]{Convention}

\begin{document}
\abovedisplayskip=6pt plus3pt minus3pt
\belowdisplayskip=6pt plus3pt minus3pt
\title[Topological obstructions to nonnegative curvature]{\bf Topological 
obstructions\\ to nonnegative curvature\rm}
\author{Igor Belegradek}
\address{Department of Mathematics and Statistics\\ McMaster University\\
Hamilton, ON, L8S 4K1, Canada}
\email{belegi@math.mcmaster.ca}
\author{Vitali Kapovitch}
\address{Department of Mathematics\\ University of Pennsylvania\\
 Philadelphia, PA 19104}
\email{vitali@math.upenn.edu}
\date{}
\begin{abstract} We find new obstructions to the existence
of complete Riemannian metric of nonnegative sectional curvature on
manifolds with infinite fundamental groups. In particular, we
construct many examples
of vector bundles whose total spaces admit no
nonnegatively curved metric.
\end{abstract}
\maketitle
\section{Introduction}
According to the soul
theorem of J.~Cheeger and D.~Gromoll a complete open manifold of
nonnegative sectional curvature is diffeomorphic to the total space of the
normal bundle of a compact totally geodesic submanifold
which is called the {\it soul}.
One of the harder questions in the subject of
is what kind of normal bundles can occur.

Cheeger and Gromoll also proved that a finite cover of any closed
nonnegatively curved manifold (throughout the paper by a nonnegatively
curved manifold we mean a complete Riemannian manifold of nonnegative {\it
sectional} curvature) is diffeomorphic to a product of a torus and
a simply-connected closed nonnegatively curved manifold.
It turns out that a similar statement holds
for open complete nonnegatively curved manifolds
(see~\cite{Wilk} and section~\ref{sec:split}
where a Ricci version of the statement is proved).

We use this fact to find
new obstructions to nonnegative curvature and build
many examples of vector bundles whose total spaces admit no complete
metric of nonnegative curvature. Basic obstructions
are provided by the following proposition (of which there
is a more general version incorporating the Euler class).

\begin{prop}\label{intro: prop polynomial in pontr}
Let $N$ be an open complete nonnegatively curved manifold
such that $Q(TN)\neq 0$ for some polynomial $Q$ in rational
Pontrjagin classes.
Then $Q(T\tilde{N})\neq 0$
for the universal (and hence any) cover $\pi\co \tilde N\to N$.
\end{prop}

Note that~\ref{intro: prop polynomial in pontr}
is true for {\it finite} covers without any curvature assumptions
(because finite covers induce injective maps on rational cohomology).
In general, the results of this paper
are only interesting for manifold with infinite
fundamental groups.

Previously, obstructions to the existence of nonnegatively curved
metrics on vector bundles were only known for a
flat soul ~\cite{OW}. No obstructions are known when the soul is
simply-connected. Examples of  nonnegatively curved metrics on
vector bundles can be found in~\cite{Che,Ri,Yan,GZ, GZ2}.

\begin{cor}\label{intro: C-times-flat}
Let $\eta$ be a vector bundle over a closed
smooth manifold $C$ and let
$\xi$ be a vector bundle over a closed flat manifold $F$
such that the total space of $\eta\times\xi$
admits a complete nonnegatively curved metric.
Then $\xi$ becomes stably trivial after passing to a finite cover.
Furthermore, if either $rank(\eta)=0$, or $\eta$ is orientable and has nonzero
rational Euler class, then $\xi$ becomes trivial in a finite cover.
\end{cor}
%

Note that a vector bundle over a flat manifold $F$ becomes
trivial in a finite cover iff its rational Euler and
Pontrjagin classes vanish. Similarly,
a bundle over $F$ is stably trivial in a finite cover iff
its rational Pontrjagin classes vanish (see~\ref{uniq}).

In case $C$ is a point~\ref{intro: C-times-flat} says that any
vector bundle $F$ with nonnegatively curved total space becomes
trivial in a finite cover. Also since the Euler and Pontrjagin
classes determine a vector bundle up to finite
ambiguity (see e.g.~\cite{Bel}), in every rank
there are only finitely many
vector bundles over $F$ with nonnegatively curved total spaces.
Thus,~\ref{intro: C-times-flat} is a 
generalization of the main result of~\cite{OW}.

To see how~\ref{intro: C-times-flat} works,
note that if $T$ is a torus of dimension $\ge 4$, then there are infinitely
many
vector bundles over $T$ of every rank $\ge 2$ with (pairwise) different
first
Pontrjagin classes. Also there are infinitely many rank $2$
vector bundles over $T^2$ and $T^3$ with different Euler classes.
We now deduce the following.

\begin{cor}\label{intro: any-base}
Let $B$ be a closed nonnegatively curved manifold. If $\pi_1(B)$
contains a free abelian subgroup of rank four (two, respectively),
then for each $k\ge 2$ (for $k=2$, respectively) there exists a
finite cover of $B$ over which there exist infinitely many rank
$k$ vector bundles  whose total spaces admit no nonnegatively
curved metrics.
\end{cor}

By contrast, any vector bundle over $S^2\times S^1$ admits a
nonnegatively curved metric as we observe in~\ref{s1s2}.
Thus~\ref{intro: any-base} cannot be
generalized to the case when $\pi_1(B)$ is virtually-$\mathbb Z$.

Passing to finite covers in~\ref{intro: C-times-flat}
and~\ref{intro: any-base} seems necessary,
in general, in order to obtain bundles without
nonnegatively curved metrics.
For example, one can easily construct
flat $SO(n)$ vector bundles over a torus
with nonzero Stiefel-Whitney classes, and obviously
their total spaces are complete flat manifolds.
Here is an example when we get a complete picture without
passing to a finite cover.

\begin{cor}\label{S1timesS3}
Let $\xi$ be a vector bundle over $S^3\times S^1$
whose total space has a nonnegatively curved metric.
Then either $\xi$ is the trivial bundle or $\xi$ is the
product of a trivial bundle over $S^3$ and the M\"obius band
bundle over $S^1$.
\end{cor}

We emphasize that our method does
not apply when $B$ is simply-connected,
or more generally if after passing to a finite cover $C\times T\to B$
the bundle $\xi$ becomes a pullback of a bundle over $C$
via the projection $C\times T\to C$.
(Here, and until the end of the section $C$ is a simply-connected
manifold and $T$ is a torus.)
For instance, if $B$ is a closed flat manifold which is an odd-dimensional
rational homology sphere~\cite{Sz}, then any vector bundle
over $B$ becomes trivial in a finite cover and it is unclear
whether there are bundles over $B$ which are not nonnegatively curved.

A reasonable goal is to find
an example of a rank $k$ vector
bundle over $C\times T$ with no nonnegatively curved metric,
whenever there is a rank $k$ vector bundle over $C\times T$ that does not
become the pullback of a bundle over $C$ in a finite cover.
This is achieved in~\ref{intro: any-base} when $\dim(T)\ge 4$.
Otherwise, the answer may depends on the topology of $C\times T$.
For example, any bundle of rank $\ge 3$ over $2$-torus
becomes trivial, and hence nonnegatively curved,
in a finite cover.

While we do not quite settle the case $\dim(T)<4$, we
get various partial results.
For instance, given a closed orientable
$2n$-manifold $B$ and an integer $d\neq 0$,
there always exists a map $f\co B\to S^{2n}$ of degree $d$.
Then, if $\pi_1(B)$ is infinite, we show that
the total space of the pullback bundle $f^\#TS^{2n}$
admits no complete metric of nonnegative curvature.
To state further results we need to review
some basic bundle theory.

By a simple obstruction-theoretic argument $H^{even}(C\times T,
C)=0$ implies that any vector bundle over $C\times T$ becomes the
pullback of a bundle over $C$ after passing to a finite cover.
This is the case, for example, for bundles over $CP^n\times S^1$.
However, once $H^{even}(C\times T, C)\neq 0$ we immediately get  a
bundle with no nonnegatively curved metric.

\begin{cor}\label{intro: euler}
If $H^{2i}(C\times T, C)\neq 0$ for some $i>0$,
then there exist infinitely many
rank $2i$ vector bundles over $C\times T$ with different Euler classes
whose total spaces are not nonnegatively curved.
\end{cor}

The Euler class is unstable and, in fact, the bundles constructed in
the proof of~\ref{intro: euler}
become pullbacks of bundles over $C$ after taking Whitney sum
with a trivial line bundle and passing to a finite cover.

To get examples that survive stabilization one has to
deal with Pontrjagin classes which live in $H^{4*}(C\times T)$.
Generally, if $H^{4*}(C\times T, C)=0$, then after adding a trivial
line bundle, any vector bundle over $C\times T$
becomes the pullback of a bundle over $C$ in a finite cover.
If $H^{4*}(C\times T, C)\neq 0$, one hopes to
find a vector bundle without nonnegatively curved metric
that survives stabilization and passing to finite covers.
We do this in several cases,
the simplest being when the rank of the bundle is $\ge\dim(C)$
(see section~\ref{sec: appl-vector} for other results involving
various assumptions on Pontrjagin classes of $TC$).

\begin{cor}
If $H^{4i}(C\times T, C)\neq 0$ for some $i>0$, then for each
$k\ge\dim(C)$ there exist infinitely many rank $k$ vector bundles
over $C\times T$ with different Pontrjagin classes
whose total spaces admit no metric of nonnegative curvature.
\end{cor}

The main geometric ingredient
of this paper is that a finite cover of any complete
nonnegatively curved manifold $N$ is diffeomorphic to a product
of a torus $T$ and a simply connected manifold $M$ and this
diffeomorphism can be chosen to take a soul $S$ to
the product of $T$ and a simply-connected submanifold of $M$.
There is also a Ricci version of this statement described in
section~\ref{sec:split}.
For example, the above conclusion holds if
$N$ has nonnegative Ricci curvature, $S$ is an isometrically
embedded compact submanifold of $N$ such that the inclusion
$S\hookrightarrow N$ induces an isomorphism of fundamental groups,
and either $S$ is totally convex, or there exists a distance nonincreasing
retraction $N\to S$.
In particular, all the theorems stated above hold in these
cases.

Our methods also yield  obstructions to existence of
metrics of nonnegative Ricci curvature on {\it closed} manifolds
(after all,~\ref{intro: prop polynomial in pontr}
can be applied to closed manifolds).
Here is an example.
It was shown in~\cite{GW} that the total space of the sphere bundle
associated with the normal bundle to the soul has a nonnegatively
curved metric. Thus, potentially, sphere bundles provide a good source of
closed nonnegatively curved manifolds. Among other things,
we prove the following.

\begin{cor}
Let $\xi$ be a bundle over a flat manifold $F$
with associated sphere bundle $S(\xi)$ and let
$C$ be a closed smooth simply-connected manifold.
If $C\times S(\xi)$ admits a metric of nonnegative Ricci curvature , 
then $\xi$ becomes trivial in a finite cover.
\end{cor}

Finally note that obstructions to to the existence of nonnegatively curved 
metrics
on total spaces of vector bundles give rise to
obstructions to the existence of
$G$-invariant nonnegatively curved metrics
on the associated $G$-principal bundles.
Indeed, any vector bundle $\xi$ with a structure group $G$
can be written as $(P\times\mathbb R^k)/G$
where $P$ is a principal $G$-bundle and $G$ acts on
$\mathbb R^k$ via a representation $G\to SO(k)$.
By the O'Neill curvature submersion formula,
if $P$ has a $G$-invariant nonnegatively curved metric,
then so does the total space of $\xi$.

The structure of the paper is as follows.
Section~\ref{sec:split} contains the above mentioned
splitting theorem for nonnegatively curved manifolds.
Section~\ref{sec: basic-obstr} summarizes
the obstructions to nonnegative curvature coming from
the splitting theorem.
In  section~\ref{sec: prodbundles} we develop
general existence and uniqueness results for bundles over
$C\times T$.
Section~\ref{sec: appl-vector} contains concrete examples of
vector bundles with no nonnegatively curved metrics.
Various obstructions to the existence of  metrics of nonnegative Ricci
curvature on sphere bundles
are described in  section~\ref{sec: spherebundles}.
Theorem~\ref{S1timesS3} is proved in the section~\ref{sec: s1s3}.

We are grateful to William Goldman, Burkhard Wilking, and Wolfgang
Ziller for many illuminating conversations. The first author is
thankful to the Geometry-Topology group of the McMaster University
for support and excellent working conditions.

\section{Splitting in a finite cover}
\label{sec:split}
Cheeger and Gromoll proved in~\cite{CG} that a finite cover of a
closed nonnegatively curved manifold is diffeomorphic to a product
of a torus and a simply connected manifold. The main geometric
tool we employ in this paper is the following generalization of
this result to open manifolds.
\begin{lem}\label{mainlemma}
Let $(N,g)$ be a complete nonnegatively curved manifold. Then
there exists a finite cover $N'$ of $N$ diffeomorphic to a product
$M\times T^k$ where $M$ is a complete open simply connected
nonnegatively curved manifold. Moreover, if $S^\prime$ is a soul of $N'$,
then this diffeomorphism can be
chosen in such a way that it takes $S'$ onto
$C\times T^k$ where $C$ is a soul of $M$.
\end{lem}
After obtaining this result we have learned that it follows from a
more general theorem which was proved earlier by B.
Wilking~\cite{Wilk}. We then realized
that our proof of~\ref{mainlemma} in fact gives  the following
stronger statement.
\begin{prop}\label{riccisplit1}
Let $(N,g)$ be a complete manifold of
nonnegative Ricci curvature.
Let $q\co \tilde{N}\to N$ be the universal
cover of $N$ and  let $\rho\co \pi\to Iso(\tilde{N})$ be the deck
transformation representation of $\pi=\pi_1(N)$.

Suppose that there exists
  a closed manifold $S\subset N$
isometrically embedded into $N$ such that the
inclusion $S\hookrightarrow N$ induces an isomorphism of the
fundamental groups, and any line
in $\tilde{S}=q^{-1}(S)$ with respect to the
induced metric from $\tilde{N}$ is also a line in $\tilde{N}$.

Then $\pi$ is  virtually  abelian and, if $\pi$ has no torsion, then
there exists a smooth path
$\rho(t)\co [0,1]\to\mathrm{Hom}(\pi,Iso(\tilde{N}))$
such that
\begin{enumerate}[(i)]
\item $\rho(0)=\rho$;
\item for each $t$ the action of $\pi$ on $\tilde{N}$ is free and
properly discontinuous;
\item \label{spl2} A finite cover of $N_1=\tilde{N}/\rho(1)(\pi)$ splits
isometrically as  $M\times T^k$ where $k=rank(\pi)$;
\item There is a family of closed submanifolds $S_t\subseteq
N_t=\tilde{N}/\rho(t)(\pi)$ such that
\begin{enumerate}[(a)]
\item $S_0=S$
\item
Under the splitting from (\ref{spl2})
the cover of $S_1$ corresponds to the Riemannian product
 $C\times T^k\subset M\times T^k$ where $C$ is a closed isometrically
 embedded submanifold of $M$.
\item for each $t$ there exists a diffeomorphism $\phi_t\co
(N_t,S_t)\to (N_0,S_0)$
\end{enumerate}
\end{enumerate}
\end{prop}
The assumption that any line in $\tilde{S}=q^{-1}(S)$
is also a line in $\tilde{N}$ is satisfied  if $S$ is
totally convex in $N$ or if there is a distance nonincreasing
retraction $N\to S$. Both of these conditions are true if $N$ is
an open manifold of nonnegative sectional curvature and
$S\subseteq N$ is its soul. In this case one can also describe the
souls of the deformed manifolds $N_t$. Namely we have the
following
\begin{prop}\label{split}
Let $(N,g)$ be a complete nonnegatively curved manifold with a
free abelian fundamental group $\pi$.
Let $\rho\co \pi\to Iso(\tilde{N})$ be the deck transformation
representation of $\pi$. Then there exists a smooth path
$\rho(t)\co [0,1]\to\mathrm{Hom}(\pi,Iso(\tilde{N}))$ such that
in addition to (i)-(iv) of~\rm\ref{riccisplit1}\it\ the following
holds.

(v) If $S$ is a soul of $N$, then there exists an isometric splitting
$\tilde{N}=M\times\mathbb{R}^k$ where $k=rank(\pi)$ and a soul
$C$ of $M$ such that, for every $t\in[0,1]$, the projection $S_t$
of $C\times\mathbb{R}^k$ to $N_t=\tilde{N}/\rho(t)(\pi)$ is a
soul of $N_t$. Also for each $t$, there exists a diffeomorphism
$\phi_t\co (N_t,S_t)\to (N_0,S_0)$.
\end{prop}
\begin{rmk} The above mentioned result of Wilking~\cite{Wilk}
implies the existence of the deformation $\rho(t)$ as
in~\ref{riccisplit1} for an arbitrary virtually abelian group. He
also gives an upper bound on the order of the covering in question
in terms of $\pi$ and the number of connected components of
$Iso(M)$. Nevertheless, we will present our proof of~\ref{riccisplit1}
for it is considerably easier than the one in~\cite{Wilk}.
(In fact, our proof is very similar to the original argument
of Cheeger and Gromoll in the closed manifold case.)
Besides, the statements of~\ref{riccisplit1} and~\ref{split}
are tailored to our applications, for example
the parts $(iv)-(v)$ are not discussed
in~\cite{Wilk}.
\end{rmk}
\begin{proof}[Proof of~\ref{riccisplit1}]
Let $q\co \tilde{N}\to N$ be the universal cover of $N$ and let
$\tilde{S}=q^{-1}(S)$. Then since inclusion $S\hookrightarrow N$
induces an isomorphism of the fundamental groups
$q|_{\tilde{S}}\co \tilde{S}\to S$ is the universal cover of $S$.

Let $\tilde{S}=C\times\mathbb{R}^k$ be the de Rham
decomposition of $\tilde{S}$ so that $C$ does not split off a
Euclidean factor. We claim that $C$ is compact.
(Indeed, suppose $C$ is not compact.
Then $C$ contains a ray $\gamma$.
Since $S$ is compact, there exists a point $p\in S$ such that
$q(\gamma(i))\to p$ as $i\to \infty$.
Let $\tilde{p}$ be a point in $q^{-1}(p)$.
By above there exists a sequence $g_i\in \pi$ such that
$g_i(\gamma(i))\to\tilde{p}$. Passing to a subsequence we can
assume that $g_i(\gamma'(i))\to v\in T_{\tilde{p}}C$.
Just as in~\cite{CG2} we readily conclude that
$\sigma(t)=\exp(tv)\co\mathbb{R}\to\tilde{S}$
is a line in $C\subset\tilde{S}$.
By assumption $\sigma(t)$ is also a line in  $\tilde{N}$.
Therefore, the splitting theorem~\cite{CG2} implies
$\tilde{N}$ splits off $\sigma(t)$ isometrically.
So $v$ is invariant under the $Hol(\tilde{N})$ and hence
under $Hol(\tilde{S})$ which contradicts the fact that $C$ does not
split off a Euclidean factor.)

Since any isometry of
$\tilde{S}$ takes lines to lines, the isometry group
$Iso(\tilde{S})$ splits as a direct product
$Iso(\tilde{S})\cong Iso(C)\times Iso(\mathbb{R}^k)$.
Therefore, the natural deck transformation action of $\pi$ on
$\tilde{S}$ gives a monomorphism
$\rho=(\rho_1,\rho_2)\co\pi\to Iso(C)\times Iso(\mathbb{R}^k)$.

Since $Iso(C)$ is compact and $\pi$ is discrete, the group
$\rho_2(\pi)$ is a discrete subgroup of $Iso(\mathbb{R}^k)$.
Also $\ker (\rho_2)$
is compact and hence it is finite. Thus, $\pi$ is an extension of a
finite group by a crystallographic one.
It is well-known (see e.g. the proof of~\cite[Thm 2.1]{Wilk} )
that any such a group is virtually abelian.

Now suppose that $\pi$ is free abelian.
Then $\rho_2(\pi)$ is a discrete torsion-free
subgroup of $Iso(\mathbb{R}^k)$, in particular, it acts
on $\mathbb{R}^k$ by translations and
$\mathbb{R}^k/\rho_2(\pi)$ is isometric to a flat torus $T^k$.

By above the splitting $\tilde{S}=C\times\mathbb{R}^k$ is just a
part of a bigger isometric splitting
$\tilde{N}=\tilde{M}\times\mathbb{R}^k$ where
$\tilde{M}$ is a complete open simply connected
manifold  containing $C$ as an isometrically
embedded submanifold.

Since the action of $\pi$ on $\tilde{N}$ leaves $\tilde{S}$ invariant, it
sends lines parallel to $\mathbb{R}^k$ into lines parallel to
$\mathbb{R}^k$.
Hence the map $\rho $ is a restriction of a natural monomorphism
$\pi\to Iso(\tilde{M})\times Iso(\mathbb{R}^k)$  which with a slight abuse
of notations we will still denote by $\rho=(\rho_1,\rho_2)$.
In fact, the image of $\rho_1$ lies in the subgroup
$G\unlhd Iso(\tilde{M})$ of isometries leaving $C$ invariant.
Since $C$ is compact it follows that $G$ is a {\it compact}
subgroup of $Iso(\tilde{M})$.

Next consider the homomorphism $\rho_1\co \pi\to G$.  Let $H$ be
the closure of $\rho_1(\pi)$ in $G$. Then $H$ is a compact abelian
subgroup of $G$. Let $H_0$ be the identity component of $H$.
Consider the short exact sequence $1\to H_0\to H\to \Gamma\to 1$
where $\Gamma=H/H_0$ is a finite abelian group. We claim that this
sequence splits and hence $H\cong H_0\times \Gamma$.

Indeed, the group $\Gamma$ is a product of finite cyclic subgroups
and, since $H$ is abelian, it is enough to define the splitting on
generators of these subgroups. Let $g\in\Gamma$ be a generator of
order $m$ and let $\bar{g}\in H$ be a preimage of $g$. The
endomorphism of $H$ sending $x$ to $x^m$, takes $\bar g$ to $H_0$,
and maps $H_0$ {\it onto} itself. Hence, there is $h\in H_0$ such
that $h^m=\bar g^m$, and we can define a splitting by mapping $g$
to $\bar g\cdot h^{-1}$.

Thus, $\rho_1\co \pi \to H\cong H_0\times\Gamma$ can be written as
a product of two representations $\rho'\co \pi \to H_0$ and
$\rho''\co \pi \to \Gamma$. Since $H_0\cong T^l$, the
representation variety $\mathrm{Hom}(\pi, H_0)$ is diffeomorphic
to a torus $T^{kl}$. Hence, we can find a smooth deformation
$\rho_1'(t) \in \mathrm{Hom}(\pi,H_0)$ such that
$\rho_1'(0)=\rho'$ and $\rho_1'(1)=1$, the trivial representation.

Crossing $\rho_1'(t)$ with $\rho''$ and $\rho_2$, we obtain a
smooth path $\rho(t)\in\mathrm{Hom}(\pi,Iso(\tilde{N}))$
such that $\rho(0)=\rho$ and $\rho(1)=1\times \rho''\times\rho_2$.
For every $t$ the action of $\pi$ on $\tilde{N}$ via
$\rho_t$ is free and properly discontinuous because so is the
action of $\pi$ on $\mathbb{R}^k$. Therefore, we get a smooth
family of  manifolds of nonnegative Ricci curvature
$N_t=\tilde{N}/\rho(t)(\pi)$ with $N_0=N$. We also
get the family $S_t=\tilde{S}/\rho(t)(\pi)\subseteq N_t$ of closed
isometrically embedded submanifolds
with $S_0=S$.

The finite cover of $N_1$ corresponding to the kernel
of $\rho''$ splits isometrically as  $M\times T^k$.
Under the splitting the cover of $S_1$ corresponds to the
Riemannian product $C\times T^k\subset M\times T^k$.

By the (relative) covering homotopy theorem the family
$(N_t,S_t)$, considered as a bundle over $[0,1]$ is smoothly
isomorphic to the trivial bundle $[0,1]\times (N,S)$. In
particular, all $S_t$'s are mutually diffeomorphic and, moreover,
have isomorphic normal bundles.
\end{proof}
\begin{proof}[Proof of~\ref{split}]
Let $S$ be a soul of $N$ and let $p\co \tilde{N}\to N$ be the
universal cover of $N$.
By the Cheeger-Gromoll soul theorem~\cite{CG}
$S$ is totally convex and the inclusion $S \hookrightarrow
N$ is a homotopy equivalence. Thus,~\ref{riccisplit1} applies and
it only remains to deduce (v).

Let $h\co N \to \mathbb{R}$ be the Cheeger-Gromoll exhaustion
function generating $S$ and let
$\tilde{h}=h\circ q \co\tilde{N}\to\mathbb{R}$
be its lift to the universal cover
$\tilde{N}$. Clearly, $\tilde{h}$ is convex.
Moreover, since every line in $\tilde{N}$ parallel to
$\mathbb{R}^k$ projects to an infinite geodesic lying in a compact
set, $\tilde{h}$ is constant along any such line.
Hence $\tilde{h}$ is given by the formula $\tilde{h}(m,t)=\bar{h}(m)$ for
some
convex function $\bar{h}\co M \to \mathbb{R}$. It is easy to see
that $\bar{h}$ is an exhaustion function. Let
$C\subseteq M$ is the soul generated by $\bar{h}$.

By construction $\bar{h}$ is invariant under the
action of $\rho_1(\pi)$ and, hence, under the action of $H$.
In particular, $\bar{h}$ is invariant under the action of
$\rho_1(t)(\pi)$ for any $t$. Therefore, $\bar{h}$
descends to a well defined
convex exhaustion function $h_t\co N_t\to \mathbb{R}$ generating
the soul $S_t=(C\times\mathbb{R}^k)/\rho(t)(\pi)$.
\end{proof}

\begin{rmk}\label{split: comp-isometry}
Actually, it follows from the proof of~\ref{riccisplit1} that some
versions of~\ref{riccisplit1} and~\ref{mainlemma} hold without any
curvature assumptions. For example, let $N$ be a complete
Riemannian manifold whose universal cover is isometric to
$M\times\mathbb{R}^n$ where $Iso(M)$ is compact.
Then a finite cover of $N$ is
diffeomorphic to the product $M\times T^k\times\mathbb R^{n-k}$.
See~\cite{Wilk} for a stronger result.
\end{rmk}

\section{Basic obstructions}
\label{sec: basic-obstr}
In this section we obtain simple topological obstructions
to nonnegative curvature coming from the results of
the section~\ref{sec:split}.

In this section and throughout the rest of this
paper we use the notation $e$ for the
Euler class, $p_i$ for the $i$th Pontrjagin class, and
$p=\sum_{i\ge 0}p_i$ for total Pontrjagin class.
Unless stated otherwise, all the characteristic classes
live in cohomology with rational coefficients.
(However, it is useful to keep in mind that $e$ and $p_i$
are in fact integral classes, that is they lie in the image
of $H^\ast(B,\mathbb Z)\to H^\ast(B,\mathbb Q)$.)

Let $S$ be a closed manifold smoothly embedded
into an open manifold $N$ such that
the inclusion $S\hookrightarrow N$ induces an isomorphism
of fundamental groups.
Let $q\co\tilde{N}\to N$ be the universal cover of $N$;
then $q\co\tilde{S}=q^{-1}(S)\to S$ is the universal cover of $S$.
Assume that after passing to a finite cover
$N$ becomes diffeomorphic to $M\times T$
where $\pi_1(M)=1$ and $T$ is a torus of positive dimension,
Further, suppose that this diffeomorphism takes (a finite cover of)
$S$ onto $C\times T$ where $C$ is a
submanifold of $M$.
Denote the normal bundles of $S$ in $N$ by $\nu_S$.

\begin{lem}\label{observ: polyn-univcover}
Suppose there is a polynomial $Q$ with rational
coefficients such that $Q(e(\nu_S),p_1(TN|_S),p_2(TN|_S),...)\neq 0$
where $\nu_S$ is assumed to be oriented if $Q$ depends on $e$.
Then $Q(e(q^\#\nu_S),p_1(T\tilde{N}|_{\tilde{S}}),
p_2(T\tilde{N}|_{\tilde{S}}),...)\neq 0$.
\end{lem}
\begin{proof}
Note that $Q_S=Q(e(\nu_S),p_1(TN|_S),p_2(TN|_S),...)\neq 0$
remains true after passing to any finite cover
because finite covers induce injective maps on
rational cohomology.
Thus, we can assume  without loss of generality
that $N$ is diffeomorphic to $M\times T$ as above and
this diffeomorphism identifies $S$ with $C\times T$.

Then the normal bundle $\nu_C^M$ of $C$ in $M$
is the pullback of $\nu_S$
via the inclusion $i_C\co C\to S$ and, also
$\nu_S$ is the pullback of $\nu_C^M$ via the projection
$\pi_C\co S\to C$.
Similarly, since $T$ is parallelizable,
$TN|_C$ is stably isomorphic to $i_C^\#TN|_S$
and  $TN|_S$ is stably isomorphic to $\pi_C^\#TN|_C$.
In particular,
$$Q_C=Q(e(\nu_C^M),p_1(TN|_C),p_2(TN|_C),...)=
i_C^\ast Q_S$$ and
$\pi_C^\ast Q_C=Q_S$. The latter implies that
$Q_C\neq 0$.

Since $C$ is simply-connected, $i_C$ factors through
the universal covering $q\co\tilde S\to S$.
In particular, $q^\ast Q_S\neq 0$ as desired.
\end{proof}

\begin{rmk}
Clearly,~\ref{observ: polyn-univcover} remains true
for any (not necessarily universal) cover $q^\prime$
in place of $q$ because
$q^{\prime\ast} Q_S=0$ implies $q^\ast Q_S=0$.
\end{rmk}

\begin{rmk}
A particular case of~\ref{observ: polyn-univcover}
remains true even without mentioning $S$.
Namely, assume only that $N$ is an open manifold
whose finite cover is diffeomorphic to $M\times T$.
Let $Q$ be a polynomial in rational
Pontrjagin classes such that $Q(TN)\neq 0$.
Then the same proof implies $Q(T\tilde{N})\neq 0$.
This applies to the geometric situation discussed
in~\ref{split: comp-isometry}.
\end{rmk}

\paragraph{Base versus soul.}
Now we specialize to the case when
$N$ is the total space
of a smooth vector bundle over a closed manifold $B$.
We identify $B$ with the zero section.
Then the universal cover of $B$ is
$q\co\tilde B=q^{-1}(B)\to B$.
Assume also that the inclusion $S\hookrightarrow N$
is a homotopy equivalence.

First, we need to see how the characteristic classes
of the normal bundles to the $S$ and $B$ are related.
The homotopy equivalence $h\co B\to S$
(defined as the composition of the inclusion
$B\hookrightarrow N$ and a homotopy
inverse of $S\hookrightarrow N$) clearly has the
property that $h^\ast p_i(TN|_S)=p_i(TN|_B)$ for any $i$.

Furthermore, if the manifolds $N$, $B$, $S$ are oriented, then
for the rational Euler class we have
$h^\ast e(\nu_S)=\deg(h)e(\nu_B)$.
(Indeed, suppressing the inclusions we have
$\langle h^\ast e(\nu_B),\alpha\rangle =
\langle e(\nu_B),h_\ast\alpha\rangle=
\langle [B], h_\ast\alpha\rangle=
\langle \deg(h)[S],\alpha\rangle=
\deg(h)\langle e(\nu_S),\alpha\rangle$.)

By possibly changing orientation on $S$ we can
arrange that $\deg(h)=1$ so that $h^\ast e(\nu_B)=e(\nu_S)$.
Thus, since $h^\ast$ is an algebra homomorphism, we get
$$h^\ast Q(e(\nu_B),p_1(TN|_{B}),p_2(TN|_{B}),...)=
Q(e(\nu_S),p_1(TN|_{S}),p_2(TN|_{S}),...).$$

\begin{prop}\label{observ: base-soul-cover}
Let $\xi$ be a vector bundle over a closed smooth manifold
$B$ whose total space $N$ admits a complete Riemannian metric
of nonnegative sectional curvature.
Suppose there is a polynomial $Q$ with rational
coefficients such that $Q(e(\xi),p_1(TN|_B),p_2(TN|_B),...)\neq 0$
where $\xi$ is assumed to be oriented if $Q$ depends on $e$.
Then $Q(e(q^\#\xi),p_1(T\tilde{N}|_{\tilde{B}}),
p_2(T\tilde{N}|_{\tilde{B}}),...)\neq 0$.
\end{prop}
\begin{proof}
Let $S$ be a soul of $N$. By above
the homotopy equivalence $h$
takes $Q(e(\xi),p_1(TN|_B),p_2(TN|_B),...)$ to
$Q(e(\nu_S),p_1(TN|_S),p_2(TN|_S),...)$ hence the
latter is nonzero.

By~\ref{observ: polyn-univcover} we have
$Q(e(q^\#\nu_S),p_1(T\tilde{N}|_{\tilde{S}}),
p_2(T\tilde{N}|_{\tilde{S}}),...)\neq 0$.
Let $\tilde h$ be the lift of $h$ to the universal covers;
note that $\tilde h$ is a homotopy equivalence.
Then by commutativity
$\tilde h^\ast e(q^\#\nu_B)=e(q^\#\nu_S)$
and $\tilde h^\ast p_i(T\tilde{N}|_{\tilde B})=p_i(T\tilde N|_{\tilde S})$.
So the homotopy inverse of $\tilde h$ takes
$Q(e(q^\#\nu_S),p_1(T\tilde{N}|_{\tilde{S}}),
p_2(T\tilde{N}|_{\tilde{S}}),...)$ to the corresponding polynomial for
$\tilde B$ which is therefore nonzero as claimed.
\end{proof}

\begin{rmk}
The statement of~\ref{observ: base-soul-cover} becomes especially
simple if $p(TB)=1$. Indeed, it implies that
$p(TN|_{B})=p(\xi\oplus TB)=p(\xi)p(TB)=p(\xi)$. We also get
$p(T\tilde B)=1$ which implies $p(T\tilde N|_{\tilde
B})=p(q^\#\xi)$.
\end{rmk}

We shall often use the following variation of~\ref{observ: base-soul-cover}.

\begin{prop}\label{observ: base-soul}
Let $\xi$ be an vector bundle over $B=C\times T$ where
$C$ is a closed connected smooth manifold
and $T$ is a torus.
Assume the total space $N$ of $\xi$ admits a complete
Riemannian metric of nonnegative sectional curvature.
Suppose there is a polynomial $Q$ with rational
coefficients such that $Q(e(\xi),p_1(TN|_B),p_2(TN|_B),...)\neq 0$
where $\xi$ is assumed to be oriented if $Q$ depends on $e$.
Then $Q(e(i_C^\#\xi),p_1(TN|_C),p_2(TN|_C),...)\neq 0$ where
$i_C\co C\to C\times T$ is the inclusion.
\end{prop}
\begin{proof}
The universal cover
$q\co\tilde B=\tilde C\times \mathbb R^k\to C\times T=B$
clearly factors through the
inclusion $i_C\co C\to C\times T$.
By~\ref{observ: base-soul-cover},
$Q(e(q^\#\xi),p_1(T\tilde{N}|_{\tilde{B}}),
p_2(T\tilde{N}|_{\tilde{B}}),...)\neq 0$, therefore,
$Q(e(i_C^\#\xi),p_1(TN|_C),p_2(TN|_C),...)$
must be nonzero.
\end{proof}

\section{Producing vector bundles}
\label{sec: prodbundles}
In this section we discuss some methods of building vector bundles.
We start from several general methods and then
concentrate on the case when the base is $C\times T$ where $C$ is a
finite connected CW-complex and $T$ is a torus of
positive dimension.

\begin{ex}
Let $B$ be a closed orientable
$2n$-manifold and let $\xi$ be a bundle over $S^{2n}$.
Since there always exists a degree one map $f\co B\to S^{2n}$, we
get a pullback bundle $f^\#\xi$. Now if $\xi$ has a nonzero
rational characteristic class (that necessarily
lives in $H^{2n}(S^{2n},\mathbb Q)$),
so does $f^\#\xi$ because $f$ induces an isomorphism on the $2n$-dimensional
cohomology. In particular, every even integer $2d$ can be realized
as the Euler number of a rank $2n$ bundle over $B$ (by taking
$\xi$ to be the pullback of $TS^{2n}$ via a self-map of $S^{2n}$
of degree $d$).
\end{ex}

\begin{ex}
Any element of $H^2(B,\mathbb Z)$ can be realized as
the Euler class of an oriented rank two bundle over $B$
(where $B$ is any paracompact space)~\cite[I.4.3.1]{Hirz}.
\end{ex}

\begin{ex}
If $B$ is a finite CW-complex of dimension $d$, then it is well-known
that a multiple
of any element of $\oplus_{i>0}H^{4i}(B,\mathbb Q)$ can be realized
as the Pontrjagin character of a vector bundle over
$B$ of rank $d$ (and hence of any rank $\ge d$).

In particular,
a multiple of any element $x\in H^{4k}(B,\mathbb Q)$
can be realized as the $k$th Pontrjagin class
of a bundle of rank $d$. (Indeed, let $X$ be the image of $x$
under the inclusion $H^{4k}(B,\mathbb Q)\to\oplus_{i>0}H^{4i}(B,\mathbb Q)$.
Realize a multiple of $X$ as the Pontrjagin character of a bundle.
Then this bundle has zero Pontrjagin classes $p_i$ for
$0<i<k$ and the $k$th Pontrjagin class is a multiple of $x$.)
\end{ex}

We now prove a uniqueness and existence theorem 
for vector bundles over $C\times T$.

\begin{thm}\label{uniq}
Let $C$  be a finite connected CW-complex and let $\xi$ and
$\eta$ be oriented rank $n$ vector bundles over $C\times T$ such that
\begin{equation}\label{restr}
\xi|_{C\times *}\cong \eta|_{C\times *}
\end{equation}
and
\begin{equation}\label{rational}
 \xi \text{ and }\eta\text{ have the same rational characteristic
 classes.}
\end{equation}
 Then there exists a finite cover
$\pi\co C\times T\to C\times T$ such that $\pi^\#\xi\cong\pi^\#\eta$.
\end{thm}
\begin{proof}

First, note that after passing to a finite cover
we can assume that $\xi$ and $\eta$ have the same integral cohomology
classes.
(Indeed, look for example at the integral $i$-th Pontrjagin class $p_i$.
By assumption $p_i(\xi)-p_i(\eta)$ is a torsion element of
$H^{4i}(C\times T ,\mathbb Z)$. By the K\"unneth formula we can write
$p_i(\xi)-p_i(\eta)=\sum_{j=0}^{4i}c^{4i-j}\otimes t^j$ where
$c^s\in H^s(C,\mathbb Z)$ and $t^{j}\in H^j(T^k,\mathbb Z)$.
Condition (1) implies
$t^0\otimes c^{4i}=p_i(\xi|_{C\times *})-p_i(\eta|_{C\times *})=0$.
Since any torsion element of the form
$\sum_{j=1}^{4i}c^{4i-j}\otimes t^j$ becomes zero
when mapped to an appropriate finite cover
$C\times T\to C\times T$ along $T$, the integral Pontrjagin classes
$p_i(\xi)$, $p_i(\eta)$ become equal
in such a finite cover.)

Let $f,g\co C\times T\to BSO(n)$ be the classifying
maps for $\xi$ and $\eta$ respectively.
Let $\g^n$ be the universal bundle over $BSO(n)$.
For each $i\le [n/2]$, we view the classes $p_i(\g^n)\in H^{4i}(BSO(n))$
as maps $p_i\co BSO(n)\to K(\mathbb Z,4i)$, and similarly if $n$ is even,
$e(\g^n)\in H^{n}(BSO(n))$ is thought of as a map
$e\co BSO(n)\to K(\mathbb Z,n)$.

Consider the combined map $c$ from $BSO(n)$ to the product of
Eilenberg-MacLane spaces given by the formula
$$c=(p_1,p_2...p_{(n-1)/2})\co BSO(n)\to X=\times_{s=1}^{(n-1)/2}
K(\mathbb Z,4s)\text{ if }n \text{ is odd}$$ and,
$$c=(e,p_1,p_2...,p_{n/2-1})\co BSO(n)\to X=K(\mathbb Z,n)\times
(\times_{s=1}^{n/2-1} K(\mathbb Z,4s)).$$
 if $n$ is even.

  It is well known that
$c$ is a rational homotopy equivalence (i.e. the homotopy fiber
$F$ of $c$ has finite homotopy groups) and the spaces $F, BSO(n),
X$ are simply-connected (see e.g.~\cite{Bel}).

By the condition (\ref{rational})
$p\circ f$ is homotopic to $p\circ g$. We shall now try to lift
this homotopy to the homotopy of $f$ and $g$.
%
%
We view the pair $(C\times T , C)$ (where we identify $C\times *$ with $C$)
as a relative CW-complex and we try to
construct the homotopy between $f$ and $g$ inductively on the
dimension of the skeletons.
%
%
By (\ref{restr})  we can assume that the homotopy is already
constructed on the zero skeleton $(C\times T , C)_0$.

Suppose that we have already constructed the homotopy on
$(C\times T , C)_{i-1}$ for $i>0$.
We want to show that after possibly passing
to a finite cover we can extend it over $(C\times T , C)_i$.
The relative obstruction
$O_i$ to the extension over the $i$-th skeleton lives in the
cohomology $H^i((C\times T , C),\pi_i(F))$.
Let $m=|\pi_i(F)|$ and $k=\dim(T)$; by assumption $k>0$.
Consider the $m^k$ cover
$\Pi\co C\times T\to C\times T$ given by the formula
$(c,z_1,...,z_k)\mapsto (c,z_1^m,...,z_k^m)$.
Notice that by the K\"unneth formula for the
pair $(C\times T, C)=(C,\emptyset)\times (T,*)$ we have
\begin{equation*}
H^i((C\times T , C),\pi_i(F))=
\sum_{j=0}^i  H^{i-j}((C,\emptyset),\pi_i(F))\otimes H^j((T,*),\pi_i(F))=
\end{equation*}
$$ =\sum_{j=1}^iH^{i-j}((C,\emptyset),\pi_i(F))\otimes H^j((T,*),\pi_i(F))$$
where the last equality is due
to the fact that $H^0((T,*),\pi_i(F))=0$.
Therefore, for any $\delta\in H^i((C\times T, C),\pi_i(F))$,
its pullback $\Pi^*(\d)$ is an $m$-th multiple of some class, and
thus is equal to zero. In particular, $\Pi^*(O_i)=0$. On the other
hand, by the naturality of obstructions $\Pi^*(O_i)$ is the
obstruction to extending the homotopy between $f\circ \Pi$ and
$g\circ \Pi$ over the relative $i$-skeleton $(C\times T , C)_i$.
\end{proof}

\begin{thm}\label{exist}
Let $\xi$ be an oriented rank $n$ vector bundle over a finite connected
CW-complex $C$ and let $i\co C \to C\times T $ be the
canonical inclusion onto $C\times *$.
Let $e' \in H^n(C\times T),
p_1' \in H^4(C\times T),\dots,
p_{[n/2]}'\in H^{4[n/2]}(C\times T)$ be a
collection of integral cohomology classes such that their
restrictions onto $C\times *$ give corresponding integral
characteristic classes of $\xi$ and, furthermore,
$e'=0$ if $n$ is odd and $p_{[n/2]}'=e'\cup e'$ if $n$ even.
Then there exists a finite cover
$\Pi\co C\times T\to C\times T $ and a rank $n$ vector bundle
$\eta$ over $C\times T$ such that the integral characteristic
classes of $\eta$ satisfy $e(\eta)=\Pi^*(e')$ and $p_i(\eta)=\Pi^*(p_i')$
for
$i=1,\dots,[n/2]$.
\end{thm}
\begin{proof}
Again, consider the universal fibration
$F\to BSO(n)\stackrel {c}{ \rightarrow} X$ where $X$ is the product of
appropriate Eilenberg-MacLane spaces.

The collection of characteristic classes
$(e',p_1',p_2',\dots, p_{[n/2]}')$ defines a natural map
$c'\co C\times T\to X$
(where we exclude $e'$ for odd $n$ and $p_{[n/2]}'$ for even $n$).
It suffices to
show that after passing to a finite cover there exists a lift of
this map to the map $f\co C\times T\to BSO(n)$.

By assumptions, we can construct the lift $f$ over
$C\times *$ by letting $f|_{C\times *}$ to be equal to a
classifying map of $\xi$.

Now we are faced with the relative lifting problem of extending
the lift from $C\times *$ to $C\times T$.
As in the proof of~\ref{uniq} we will proceed by induction on the dimension
of the
relative skeleton $(C\times T, C\times *)_i$.

Suppose $f$ is already defined on $(C\times T, C\times *)_{i-1}$
for some $i>0$. As before the primary obstruction $O_i$ to extending
the lift over $(C\times T, C\times *)_i$ lives in the group
$H^i((C\times T, C),\pi_{i-1}(F))$. Arguing exactly as in the
proof of~\ref{uniq}, we see that the cover
$\Pi\co C\times T\to C\times T$ given by the formula
$(c,z_1,...,z_k)\mapsto (c,z_1^m,...,z_k^m)$
where $m=|\pi_{i-1}(F)|$ has the property that
$\Pi^*(O_i)=0$. Therefore, the lift of $c'\circ \Pi$ given by
 $f\circ\Pi$ can be extended over the relative $i$-th skeleton
 of $(C\times T, C\times *)$. This completes the proof of the
 induction step and hence the proof of the theorem.
\end{proof}

\begin{rmk}
Note that by construction the cover $\Pi$ depends only on $n$ and
$\dim(C\times T)$.
\end{rmk}

\begin{rmk}
The proof of~\ref{exist} shows how to compute the characteristic classes
of a bundle with the classifying map $f$.
For example, represent $e'$ as $\sum_{j=0}^k e'_j$ where
$e'_j\in H^{n-j}(C,\mathbb Z)\otimes H^j((T,*),\mathbb Z)$. Then
the Euler class of $f$ is given by $\sum_{j=0}^k m^je'_j$ where $\dim(T)=k$.
The same result is of course true for any Pontrjagin class of $f$.

In particular, if $e'=e'_j$ for some $j>0$, then an integer
multiple of $e'$ is realized as the Euler class of some bundle
over $C\times T$.
\end{rmk}

\begin{ex}\label{exist1}
We shall often use~\ref{exist} in the following situation.
Assume $H^{4i}(C\times T , C,\mathbb Z)$ is infinite and let
$p_i^\prime\in H^{4i}((C\times T, C),\mathbb Z)$ be a nontorsion class
and $j$ be any nonzero integer.
Let $\xi$ be the trivial bundle of some rank $>2i$. Then by~\ref{exist}
there exists a bundle $\eta_j$ over $C\times T $
and a finite cover $\Pi\co C\times T\to C\times T$ such that
the restriction of $\eta_j$ to $C\times *$ is isomorphic to $\xi$
and  $\Pi^\ast (j p^\prime_i)=p_i(\eta_j)$.
Clearly, the bundles $\eta_j$ are pairwise nonisomorphic.
\end{ex}

\section{Vector bundles with no nonnegatively curved metrics}
\label{sec: appl-vector}
In this section we obtain concrete
examples of bundles without nonnegatively curved metrics.
Throughout this section $T$ is a torus of positive dimension and
$C$ is a closed connected smooth manifold.
(Note that $C$ is not assumed to be simply-connected
so the results of this section are slightly more general
than the ones stated in the introduction.)
We shall often use that the tangent bundle
to $C\times T$ is stably isomorphic to
the pullback of $TC$ via the projection $C\times T\to C$.
All (co)homology groups
and characteristic classes in this section have rational
coefficients.

\begin{cor}\label{torus-finitecov}
Let $\eta$ be a vector bundle over $C$ and let
$\xi$ be a vector bundle over $T$
such that the total space of $\eta\times\xi$
admits a complete nonnegatively curved metric.
Then $\xi$ becomes stably trivial in a finite cover.
Furthermore, if either $rank(\eta)=0$, or $\eta$ is orientable with 
$e(\eta)\neq 0$, then $\xi$ becomes trivial in a finite cover.
\end{cor}
\begin{proof}
Denote $\eta\oplus TC$ by $\eta^\prime$ so that
the tangent bundle to the total space of $\eta\times\xi$
restricted to the zero section is stably isomorphic to
$\eta^\prime\times\xi$.
Let $i$ be the largest nonnegative integer such that
$p_i(\eta^\prime)\neq 0$.
 Arguing by contradiction assume that
$p_k(\xi)\neq 0$ for some $k>0$. Using the product formula
$p(\eta^\prime\times\xi)=p(\eta^\prime)\times p(\xi)$,
we conclude that the
component of $p_{i+k}(\eta^\prime\times\xi)$ in the group
$H^{4i}(C)\otimes H^{4k}(T)$ is equal to
$p_i(\eta^\prime)\times p_k(\xi)$. Since the cross product
of nonzero classes is nonzero,
$p_{i+k}(\eta^\prime\times\xi)$ is nonzero. On the other hand, the
component of $p_{i+k}(\eta^\prime\times\xi)$ in
$H^{4i+4k}(C)\otimes H^{0}(T)$ is
$p_{i+k}(\eta^\prime)\times 1=0\times 1=0$. We now
apply~\ref{observ: base-soul} for $Q=p_{i+k}$ to get a
contradiction. By~\ref{uniq}, $\xi$ becomes stably trivial
in a finite cover.

Now assume that either $rank(\eta)=0$, or $e(\eta)\neq 0$.
By~\ref{uniq} it suffices to show that $e(\xi)=0$.
Of course, we can assume that $rank(\xi)>0$.
The pullback of $\eta\times\xi$ to $C$ has zero Euler class
because the pullback is the Whitney sum
of $\eta$ and a trivial bundle of the same rank as $\xi$.
Hence according to~\ref{observ: base-soul}, $e(\eta\times\xi)=0$.
Thus if $rank(\eta)=0$, we get $e(\xi)=0$.
Otherwise, note that $e(\eta\times\xi)=e(\eta)\times e(\xi)$ and 
since $e(\eta)\neq 0$ it implies $e(\xi)=0$ as wanted.
\end{proof}
\begin{rmk}
The assumption that $e(\eta)\neq 0$ is certainly necessary, in general.
For example, let $\eta$ be the trivial line bundle over $C$ and
let $\xi$ be the bundle over $T^{2n}$ which is the pullback of $TS^{2n}$
via a degree one map $T^{2n}\to S^{2n}$.
Then the bundle $\eta\times \xi$ is trivial so its total space is 
nonnegatively curved whenever $sec(C)\ge 0$.
\end{rmk}
%
%
\begin{thm}\label{ex: polynomialQ}
Let $H^{4i}(C\times T, C)\neq 0$ for some $i>0$ and let
$Q^\prime$ be a polynomial in rational Pontrjagin classes $p_j$ where
$0<j<i$ such that the projection of $p_i(TC)+Q^\prime (TC)$
to $H^{4i}(C)$ is zero.
Then, in each rank $>2i$, there exist infinitely many
vector bundles over $B=C\times T$ whose
total spaces are not nonnegatively curved.
\end{thm}
\begin{proof} Set $Q=p_i+Q^\prime$.
Since $H^{4i}(C\times T, C)\neq 0$, we can use~\ref{exist} to find
a vector bundle $\xi$ over $C\times T$
of rank $2i+1$ such that $p_j(\xi)=0$ for $0<j<i$ and
$p_i(\xi)$ is a nonzero class whose projection to $C$ is zero.
Using the Whitney sum formula
for Pontrjagin classes we get
$p_i(TB\oplus\xi)=p_i(TB)+p_i(\xi)$ and
$p_j(TB\oplus\xi)=p_j(TB)$ for $0<j<i$.
Thus, looking at the projection to $H^{4i}(C\times T)$,
we get $$Q(TB\oplus\xi)=p_i(TB\oplus\xi)+
Q^\prime(TB\oplus\xi)=p_i(TB)+p_i(\xi)+Q^\prime(TB)=p_i(\xi)$$
where the last equality is true because
$$p_i(TB)+Q^\prime(TB)=(p_i(TC)+Q^\prime (TC))\times 1=0\times 1=0.$$
Thus, we can apply~\ref{observ: base-soul}.
\end{proof}

\begin{cor}
Let $H^{4i}(C\times T, C)\neq 0$ for some $i>0$ and let
$ph_i(TC)=0$, where $ph_i$ is the component of the Pontrjagin character
that lives in the $4i$th cohomology.
Then, in each rank $>2i$, there exist infinitely many
vector bundles over $B=C\times T$ whose
total spaces are not nonnegatively curved.
\end{cor}
\begin{proof}
Take $Q^\prime=ph_i-p_i$.
\end{proof}

\begin{cor}\label{Q=0}
Let $H^{4i}(C\times T, C)\neq 0$ and $p_i(TC)=0$ for some $i>0$, then
there exists a vector bundle $\xi$ over $C\times T$
of any rank $> 2i$ such that
$E(\xi)$ has no metric of nonnegative curvature.
\end{cor}
\begin{proof}
Take $Q=0$.
\end{proof}

\begin{rmk}
In particular,~\ref{Q=0} shows that
if $H^{4}(C\times T, C)\neq 0$ and $p_1(TC)=0$, then
there exist infinitely many bundles of every rank $>2$ whose total
spaces are not nonnegatively curved.
\end{rmk}

\begin{cor}
If $H^{2i}(C\times T, C)\neq 0$ for some $i>0$,
then there exist infinitely many
rank $2i$ vector bundles over $C\times T$ with different Euler classes
whose total spaces are not nonnegatively curved.
\end{cor}
\begin{proof}
It follows from~\ref{exist} that there exists
a bundle $\xi$ (and, in fact, infinitely many such bundles)
of rank $2i$ over $C\times T$ such that $e(\xi)$  is a nonzero class whose
$H^{2i}(C)$ component is zero.
By~\ref{observ: base-soul} $E(\xi)$ is not nonnegatively curved.
\end{proof}

\begin{cor}\label{ex: pdual}
Let $\dim(C)=4m+2$ and $p_m(TC)\neq 0$, and let $\dim(T)\ge 2$.
Then there exist infinitely many bundles of each rank $\ge 2$
over $C\times T$ whose total spaces are not nonnegatively curved.
\end{cor}
\begin{proof}
By Poincar\'e duality find $y\in H^2(C)$ such that
$p_m(TC)y\neq 0\in H^{4m+2}(C)$. Take any $t\in H^2(T)$ and
realize $y\otimes 1+1\otimes t$ as the Euler class of an oriented rank $2$
bundle
$\xi$ over $C\times T=B$.
Then
$p_1(\xi)=(y\otimes 1+1\otimes t)^2=y^2\otimes 1+2y\otimes t+ 1\otimes t^2$.
Note that $p_m(TB)p_1(\xi)\neq 0$. (Indeed, it suffices to show that
the projection
of $p_m(TB)p_1(\xi)$ to $H^{4m+2}(C)\otimes H^2(T)$ is nonzero.
which is true because the projection is equal to $2 p_m(TC)y\otimes t$.)
Also the Whitney sum formula implies that
$p_{m+1}(TB\oplus \xi)=p_m(TB)p_1(\xi)\neq 0$ and the projection of
$p_{m+1}(TB\oplus \xi)$ to $C$ vanishes because $H^{4m+4}(C)=0$.
Hence, we are done by~\ref{observ: base-soul}.
By adding trivial bundles to $\xi$
one can make its rank arbitrary large.
\end{proof}

\begin{cor}\label{ex:add-norm-bundle}
If $H^{4i}(C\times T, C)\neq 0$ for some $i>0$, then there exist infinitely
many
vector bundles over $B=C\times T$ of any rank $\ge\dim(C)$ whose
total spaces are not nonnegatively curved.
\end{cor}
\begin{proof}
Let $\nu(C)$
be a rank $\dim(C)$ bundle over $C$ which is stably isomorphic to
stable normal bundle of $C$. By~\ref{exist} we can find a bundle $\xi$
(and, in fact, infinitely many such bundles) of rank $\dim(C)$
over $C\times T$ such that the pullback of $\xi$ to $C$ is isomorphic to
$\nu(C)$ and $p_i(\xi)$ has a nonzero projection to $H^{4i}(C\times T, C)$.

Look at the bundle $TE(\xi|_{B})\cong\xi\oplus TB$.
Note that the pullback of $\xi\oplus TB$ to $C$ is isomorphic
to $\nu(C)\oplus TC$ which is stably trivial, hence
$p_i(TE(\xi|_{C}))=0$. On the other hand the projection
of  $p(\xi\oplus TB)=p(\xi)p(TB)$ to $H^{4i}(C\times T, C)$ is equal to
$p_i(\xi)$, in particular the projection
of  $p_i(\xi\oplus TB)$ to  $H^{4i}(C\times T, C)$ is nontrivial.
By~\ref{observ: base-soul} $E(\xi)$ is not nonnegatively curved.
\end{proof}

\begin{rmk}
The method of~\ref{ex:add-norm-bundle} can be used with
some other bundles in place of $\nu(C)$.
To illustrate the idea we discuss the case when $C$ is
the total space $S(\eta)$ of the sphere
bundle associated with a vector bundle $\eta$ over $S^4$.

First, let us handle the easier case  when $e(\eta)\neq 0$.
Then $rank(\eta)$ is necessarily $4$ and it follows from the
Gysin sequence that $S(\eta)$ is a rational homology
$7$-sphere. In particular, all the Pontrjagin classes of $S(\eta)$ vanish.
Now~\ref{exist} implies that there are infinitely many rank $4$ bundles
over $S(\eta)\times T$ with nonzero $p_2$.
By~\ref{Q=0} their total spaces admit no nonnegatively curved metrics
and as usual we can add trivial bundles to make the rank $\ge 4$.

Now assume $e(\eta)=0$.
Recall that $p_1(\eta)[S^4]$ is necessarily even and
furthermore, any even integer can be realized as $p_1(\xi)[S^4]$
where $\xi$ is a $4$-bundle~\cite{Miln}.
Thus we can find a $4$-bundle $\xi'$ over $S^4$ with
$p_1(\xi')[S^4]=-p_1(\eta)[S^4]$ so that $p_1(\eta\oplus\xi')=0$.
Note that $TS(\eta)$ is stably isomorphic to the pullback of
$\eta$ via the bundle projection $\pi\co S(\eta)\to S^4$.
Setting $\xi=\pi^\#\xi'$, we get that $p_1(TS(\eta)\oplus\xi)=0$.
Since $e(\eta)=0$, the Gysin sequence implies that
$\pi$ induces an isomorphism on $H^4$, and by the Poincar\'e
duality $H^3(S(\eta))\neq 0$. Hence by~\ref{exist}
there are infinitely many rank $4$ bundles
over $S(\eta)\times T$ with nonzero $p_1$ such that
their pullback to $S(\eta)$ is $\xi$.
So the proof of~\ref{ex:add-norm-bundle} applies and
these bundles admit no nonnegatively curved metrics.

It is interesting to see whether $rank(\xi)$ can be lowered to $3$.
Recall that an integer $k$ can be realized as
$p_1(\xi)[S^4]$ for a $3$-bundle
over $S^4$ iff $k$ is a multiple of $4$~\cite{Miln}.
Thus, if $p_1(\eta)[S^4]$ is divisible by $4$, the argument
of the previous paragraph applies and we get infinitely many
$3$-bundles over $S(\eta)\times T$
with no nonnegatively curved metrics.

Now assume that $p_1(\eta)[S^4]\equiv 2 \mathrm{mod}(4)$.
We are looking for a $3$-bundle
$\xi$ over $S(\eta)$ such that $p_1(TS(\eta)\oplus\xi)=0$.
Since $e(\eta)=0$, the bundle $S(\eta)\to S^4$ has a section $s$.
Setting $\xi^\prime=s^\#\xi$, we would get a $3$-bundle $\xi^\prime$
over $S^4$ with $p_1(\eta\oplus\xi^\prime)=0$.
In particular, $p_1(\xi')[S^4]\equiv 2 \mathrm{mod}(4)$
which is impossible for a $3$-bundle.

Thus, the methods of this paper fail here. For instance,
we do not have examples of $3$-bundles over $S(\eta)\times S^1$
that admit no nonnegatively curved metrics whenever
$p_1(\eta)[S^4]\equiv 2 \mathrm{mod}(4)$.
Note that~\ref{exist} produces many  $3$-bundles over
$S(\eta)\times S^1$ which do not become pullbacks of
bundles over $S(\eta)$.
\end{rmk}

\paragraph{Metastable range, sphere bundles, and surgery.}
We now describe yet another variation of~\ref{ex:add-norm-bundle}.
When the method works, it gives a result similar to~\ref{ex:add-norm-bundle}
with sometimes lower rank.
We showed above that, under certain assumptions on the
Pontrjagin classes of $B=C\times T$, there are vector
bundles over $B$ whose total spaces admit no nonnegatively
curved metric.
Now the idea is to replace $C$ by a homotopy equivalent closed
manifold $C^\prime$ with ``nicer'' (e.g. trivial) Pontrjagin classes.
Then theorems of this section can be used to
produce a vector bundle over $C^\prime\times T$
whose total space admits no nonnegatively curved metric,
and can often use it to get a similar bundle over $C\times T$.

Indeed, let $f\co B\to B^\prime$ be a homotopy equivalence of
closed smooth manifolds and let $\xi$ be a vector bundle over
$B^\prime$ with total space $E(\xi)$.
Assume now that $2\mathrm{rank}(\xi)\ge \dim(B)+3\ge 5$, that
is, we are in the metastable range. By~\cite{Hae},
the homotopy equivalence $f\co B\to E(\xi)$ is homotopic to a smooth
embedding $e\co B\to E(\xi)$.
The above inequality implies that $\mathrm{rank}(\xi)\ge 3$,
hence by~\cite[Thm 2.2]{Sieb} $E(\xi)$ is diffeomorphic to
the total space of the normal bundle to $E(\nu_e)$.
Clearly, $E(\xi)$ is nonnegatively curved iff so is $E(\nu_e)$.

\begin{thm}
Let $T$ be a torus and let $C$ be a closed smooth manifold homotopy
equivalent to a closed manifold $C'$ such that $T$ and $C'$ satisfy the
assumptions of\rm~\ref{ex: polynomialQ}\it\ or\rm~\ref{ex: pdual}.\it\
Then, in each rank $>1+\dim(C\times T)/2$, there exist
infinitely many
vector bundles over $B=C\times T$ whose
total spaces are not nonnegatively curved.
\end{thm}
\begin{proof}
By~\ref{ex: polynomialQ} or~\ref{ex: pdual} we can find
a bundle $\xi$ over $C^\prime\times T$
whose total space $E(\xi)$ is not nonnegatively curved in any rank $>2i$.
Assume now that the rank is $>1+\dim(C\times T)/2$ (note that
$1+\dim(C\times T)/2>2i$ because $H^{4i}(C\times T, C)\neq 0$.
This puts us in the in metastable range so the
homotopy equivalence
$f\times\mathrm{id}\co C\times T\to C^\prime\times T\hookrightarrow E(\xi)$
is homotopic to an embedding whose normal bundle has total space
diffeomorphic to $E(\xi)$. Of course, the total space of this normal bundle
is not nonnegatively curved. By varying $\xi$
(or, rather, the  Pontrjagin class of $\xi$), we get infinitely
many such examples.
\end{proof}

One way to replace $C$ by a manifold $C^\prime$ with ``nicer''
Pontrjagin classes is by surgery. Namely, assume $\pi_1(C)=1$ and
let $\tau$ be a vector bundle
over $C$ so that $\tau$ and $TC$ are stably fiber homotopy equivalent.
Then, if the surgery obstruction vanishes (which always happens
if $\dim(C)$ is odd~\cite[II.3.1]{Br}), then there is a closed smooth
manifold
$C^\prime$ and a homotopy equivalence $f\co C^\prime\to C$
such that $TC^\prime$ is stably isomorphic to $f^\#\tau$.

In general, it is not easy to decide when a given
vector bundle, such as $TC$, is stably fiber homotopy equivalent
to a bundle with ``nicer'' Pontrjagin classes.
However, each bundle with ``nice'' Pontrjagin classes is usually
stably fiber homotopy equivalent to infinitely many different
bundles.

Indeed, recall that two vector bundles are stably fiber homotopy equivalent
if the corresponding spherical fibrations are stably equivalent.
The stable equivalence spherical fibrations over $C$
(or any finite simply connected cell complex) are
in one-to-one correspondence with $[C, BSG]$.
The Whitney sum gives $BSG$ and $BSO$ an $H$-group structure, and
the natural map $BSO\to BSG$ that assigns to a vector bundle
the corresponding spherical fibration induces a group homomorphism
$[C, BSO]\to [C, BSG]$. After tensoring with rational
the group $[C, BSO]$ becomes $\oplus_{i>0}H^{4i}(C,\mathbb Q)$ while
$[C, BSG]$ becomes the trivial group.
In particular, if $\oplus_{i>0}H^{4i}(C,\mathbb Q)\neq 0$,
each stable fiber homotopy equivalence class contains
infinitely many vector bundles in any rank $\ge\dim(C)$.

For example, let $C$ be the total
space of a sphere bundle over closed simply-connected
manifold $V$ associated with a vector bundle $\eta$.
(Due to~\cite{GW} such manifolds could be
a good source of nonnegatively curved manifolds.)
The bundle $TC$ is the pullback of $TV\oplus\eta$ via
the bundle projection $C\to V$, hence
$TC$ is stably fiber homotopy trivial whenever so are $TV$ and $\eta$.
This construction gives many manifolds with stably fiber homotopy trivial
tangent bundles.
%

\section{Sphere bundles with no metric of nonnegative curvature}
\label{sec: spherebundles}
It was shown in~\cite{GW} that the total space of the sphere bundle
associated with the normal bundle to the soul has a nonnegatively
curved metric. Thus, potentially, sphere bundles provide a good source of
closed nonnegatively curved manifolds.

\begin{thm}\label{sphere-euler} For $k>0$,
let $E\to F$ be a $k$-sphere Serre fibration
over a flat manifold $F$ with nonzero rational Euler class.
Let $P$ be closed smooth manifold such that
there is a map $P\to E$ that induces an isomorphism of
fundamental groups.
Then $P$ admits no metric of nonnegative Ricci curvature.
\end{thm}
\begin{proof}
Arguing by contradiction, assume that $P$ admits a metric of nonnegative Ricci 
curvature.
Pass to a finite cover $\tilde P\to P$ so that
$\tilde P$ is diffeomorphic to $C\times T$ where $C$ is
simply connected and $T$ is a torus.

Look at the corresponding covers $\tilde E\to E$ and
$\tilde F\to F$. Note that $\pi_1(\tilde F)$ is free abelian
because $\pi_1(\tilde F)$ is a torsion free group
which is the image of
a finitely generated abelian group $\pi_1(\tilde E)\cong\pi_1(\tilde P)$.
The $k$-sphere fibration $\tilde E\to \tilde F$ still has
nonzero rational Euler class since since it is a pullback of
$E\to F$ and since finite covers induce injective maps
on rational cohomology.

First consider the case $k=1$. The  circle fibration
$\tilde E\to \tilde F$ induces an epimorphism
$\phi\co\pi_1(\tilde E)\to\pi_1(\tilde F)$
of finitely generated free abelian groups. Therefore,
$\phi$ has a section. Since $\tilde E$ is aspherical,
this section is induced by a continuous map
$\tilde F\to \tilde E$ which defines a homotopy
section of the  circle fibration $\tilde E\to \tilde F$.
Thus, the Euler class must be zero which is a contradiction.

Now assume that $k>1$ so that $\tilde E\to \tilde F$
induces a $\pi_1$-isomorphism.
Then the inclusion
$T\to (C\times T)=\tilde P$ followed by the map
$\tilde P\to \tilde E\to \tilde F$ induces a $\pi_1$-isomorphism hence
is a homotopy equivalence. Let $s$ be its homotopy inverse.
Then $s$ followed by the inclusion $T\to \tilde P$ and the map
$\tilde P\to \tilde E$ is a homotopy section of the fibration
$\tilde E\to \tilde F$. The Euler class then must be zero which
gives a contradiction.
\end{proof}

\begin{rmk}
The above argument is a special case of the following phenomenon.
Suppose we have a Serre fibration $C\to P\to F$ where $F$ is a flat manifold
and $C$ is connected and simply-connected.
Look at the spectral sequence of this fibration with rational coefficients.
Then if there exists a nonzero differential, then
$P$ does not admit a nonnegatively curved metric.

Indeed, if $P$ is nonnegatively curved,
then a finite cover $\tilde P$ of $P$
splits topologically as $M\times T$ where $M$ is simply connected and
$T$ is a torus.
By naturality we can see that spectral sequence
of the pullback fibration $C\to\tilde P\to T$ also
has a nonzero differential.
Since the universal cover of $\tilde P$ is homotopy
equivalent to both $M$ and $C$, they are homotopy
equivalent to each other.
In particular $\dim H^*(M)=\dim H^*(C)$ and hence
$\dim H^*(\tilde P)=\dim H^*(M)\otimes H^*(T)=\dim (H^*(C)\otimes H^*(T))$.
On the other hand if there is a nonzero differential we should have that
$\dim H^*(\tilde P)<\dim (H^*(C)\otimes H^*(T))$ which is a contradiction.
\end{rmk}

\begin{thm}\label{sphere-pontr}
Let $E(\xi)$ be the total space  of a vector bundle
$\xi$ over a closed smooth manifold $B$ and let
$S(\xi)\to B$ be the associated sphere bundle.
Assume that $\xi$ has zero rational Euler class
and there exists a polynomial $Q$ in rational
Pontrjagin classes such that $Q(TE(\xi))\neq 0$
and $Q(T\tilde{E}(\xi))= 0$
for the universal cover $\pi\co \tilde{E}(\xi)\to E(\xi)$.
Then $S(\xi)$ admits no metric of nonnegative Ricci curvature.
\end{thm}

\begin{proof} First, we introduce several notations.
Let $q\co \tilde B\to B$ be the universal covering
and $j\co S(\xi)\to E(\xi)$ be the inclusion.
Also denote by $\pi$ and $i$ the bundle projection and the
zero section of $\xi$, respectively.

Since $S(\xi)\subset E(\xi)$ is an oriented codimension one
hypersurface, $TS(\xi)$ is stably isomorphic to $j^\#TE(\xi)$,
hence $Q(TS(\xi))=j^\ast Q(TE(\xi))$.
Also $TE(\xi)$ is isomorphic to $(i\circ\pi)^\#TE(\xi)$ since
$i\circ\pi$ is homotopic to the identity of $TE(\xi)$.
We get $Q(TS(\xi))=j^\ast\pi^\ast Q(i^\#TE(\xi))$.
By assumption $Q(TE(\xi))$, and hence $Q(i^\#TE(\xi))$ is nonzero.
Also $j^\ast\pi^\ast=(\pi\circ j)^\ast$ where $\pi\circ j\co S(\xi)\to B$
is the bundle projection. It follows from the Gysin sequence
that $\pi\circ j$ is injective in cohomology
because the kernel of $(\pi\circ j)^\ast$ consists of
the cup-multiples of the Euler class which is zero by assumption.
Thus, $Q(TS(\xi))\neq 0$.

On the other hand, the inclusion $S(q^\#\xi)\hookrightarrow E(q^\#\xi)$
takes $Q(TE(q^\#\xi))=0$ to $Q(TS(q^\#\xi))$. Hence,  $Q(TS(q^\#\xi))=0$
and we are in position to apply the theorem~\ref{observ: base-soul-cover}.
\end{proof}

\begin{cor}
Let $\xi$ be a bundle over a flat manifold $F$
with associated sphere bundle $S(\xi)$ and let
$C$ be a closed smooth simply-connected manifold.
If $C\times S(\xi)$ admits a metric of nonnegative Ricci
curvature, then $\xi$ becomes trivial in a finite cover.
\end{cor}
\begin{proof}
By~\ref{uniq} it suffices to show that $e(\xi)=0$ and $p(\xi)=1$.
Vanishing of $e(\xi)$ follows from~\ref{sphere-euler}.
Vanishing of all Pontrjagin classes
follows exactly as in
the proof of~\ref{torus-finitecov} where instead of referring
to~\ref{observ: base-soul-cover}
we use~\ref{sphere-pontr}.
\end{proof}

\section{The classification of nonnegatively curved vector
bundles over $S^1\times S^3$}
\label{sec: s1s3}
In this section we prove the theorem~\ref{S1timesS3}. Note that
the converse of~\ref{S1timesS3} is trivially true, i.e.
both the trivial bundle and the product of the trivial bundle
over $S^3$ and M{\"o}bius band
line bundle over $S^1$ are nonnegatively curved.
%
%
\begin{proof}[Proof of~\ref{S1timesS3}]
Any vector bundle over
a $4$-complex is the Whitney sum of a trivial bundle and a
bundle of rank $\le 4$, hence it
suffices to consider the bundles of
rank $k$ at most $4$.

First, assume that $\xi$ is orientable.
Let $q\co\mathbb{R}\times S^3\to S^1\times S^3$ be the universal cover of
$S^1\times S^3$. Then since any vector bundle over $S^3$ is
trivial we have that $q^\#(\xi)$ is trivial.
In particular, $p_1(q^\#(\xi))=e(q^\#(\xi))=0$.
Therefore, according to~\ref{observ: base-soul-cover}.
the classes $p_1(\xi)$ and $e(\xi)$ vanish. Thus, it suffices
to prove the following.
%
%
\begin{lem}\label{orient}
Let $\eta$ be an orientable vector bundle over $S^1\times S^3$
such that $p_1(\eta)=e(\eta)=0$. Then $\eta$ is trivial.
\end{lem}
\begin{proof}
Since
$H^1(S^1\times S^3,\mathbb Z/2\mathbb Z)=0=H^2(S^1\times S^3,\mathbb Z)$,
any rank one or rank two orientable
bundle over $S^1\times S^3$ is trivial.

Assume that $\eta$ is an orientable bundle of rank $4$.
Let $f\co S^1\times S^3\to BSO(4)$ denote a classifying map for
$\eta$, i.e $\eta\cong f^* \gamma^4$ where $\g^4$ is the universal
$4$-bundle over $BSO(4)$.
The first four homotopy groups of $BSO(4)$ are as follows: $\pi_0(BSO(4))=0,
\pi_1(BSO(4))=0, \pi_2(BSO(4))=\mathbb Z/2\mathbb Z, \pi_3(BSO(4))=0$ and
$\pi_4(BSO(4))=\mathbb Z\oplus \mathbb Z$. Consider the standard product
cell
decomposition of $S^1\times S^3$ coming from canonical cell
decompositions $S^1=e^0\bigcup e^1$ and $S^3=e^0\bigcup e^3$. Then
the $3$-skeleton of $S^1\times S^3$ is the wedge
$S^1\bigvee S^3$.
Since any orientable vector bundle over $S^1$ or $S^3$ is trivial,
$f|_{S^1\bigvee S^3}$
is homotopic to a point and therefore by the homotopy extension
property we can assume that $f$ send $S^1\bigvee
S^3$ to a point to begin with.
In other words, $f$ can be written
as a composition $f=\bar{f}\circ \pi$ where $\pi$ is the
factorization map $\pi\co S^1\times S^3\to S^1\times S^3\slash
(S^1\bigvee S^3)\cong S^4$. Since $\pi$ induces an isomorphism on
$H^4$, the bundle $\bar{f}^*(\g^4)$ has zero Euler and Pontrjagin
classes. It is a well known that a bundle over $S^4$ with
zero Euler and Pontrjagin classes is trivial.
(Indeed, the map $(e,p_1)\co\pi_4(BSO(4))\to \mathbb Z\oplus \mathbb Z$
which associates to a $4$-bundle over $S^4$ its Euler and
Pontrjagin classes is a rational homotopy equivalence. Then
the induced map on $\pi_4$ has finite, and hence trivial, kernel
because $\pi_4(BSO(4))\cong\mathbb Z\times\mathbb Z$.) Thus
$\bar{f}$, and hence $f$, is nullhomotopic.

A very similar argument shows
that any orientable $3$-bundle over $S^1\times S^3$ with zero
first Pontrjagin class is trivial. Again, everything can be reduced
to 3-bundles over $S^4$ with zero $p_1$.
The map $p_1\co BSO(3)\to K(\mathbb Z,4)$ is a rational homotopy
equivalence,
in particular, the induced map on $\pi_4$ has finite, and hence trivial,
kernel
because $\pi_4(BSO(3))\cong\mathbb Z$. Hence, only the trivial $3$-bundle
over
$S^4$
has zero $p_1$.
%
\end{proof}

Now suppose that $\xi$ is not orientable and its total
space admits a metric of nonnegative curvature.
Since $\mathbb{Z}$ has a unique subgroup of index $2$ the orientation
double cover for $\xi$ is given by the map $\pi_S=(z\to
z^2)\times \mathrm{id}\co S^1\times S^3\to S^1\times S^3$. Then the
pullback $\pi_S^\#(\xi)$ is orientable and also admits a metric of
nonnegative curvature.
By above, the pullback  $\pi_S^\#(\xi)$ is trivial.
The following lemma completes
the proof of~\ref{S1timesS3} in the nonorientable case.
\end{proof}
\begin{lem}\label{nonor}
Let $\eta$ be a nonorientable rank $k$ bundle over $S^1\times S^3$
whose orientation lift is a trivial bundle. Then $\eta$ is
isomorphic to the product $\mu^1\times\e^{k-1}$
of the M{\"o}bius band line
bundle $\mu^1$ over $S^1$ and a trivial rank $(k-1)$-bundle $\e^{k-1}$
over $S^3$.
\end{lem}
\begin{proof}
Since $H^1(S^1\times S^3,\mathbb Z/2\mathbb Z)=\mathbb Z/2\mathbb Z$,
there is only one nonorientable line bundle over $S^1\times S^3$,
namely, $\mu^1\times\e^{0}$.

\bf{Case of rank four.}\rm\
Let $f\co S^1\times S^3\to BO(4)$ be the classifying map for $\eta$ and
$f_0\co S^1\times S^3\to BO(4)$ be the classifying map for
$\mu^1\times\e^3$. We want to show that these maps are
homotopic. The same argument as
in the proof of~\ref{orient} shows that $f$ and $f_0$
are homotopic on the $3$-skeleton. Let us show that this homotopy
can be extended over the $4$-cell.

Let $\pi_B \co BSO(4) \to BO(4)$ be the canonical double cover.
Then each of the maps $f\circ\pi_S$ and
$f_0\circ \pi_S$ lifts to a map
$\tilde{f}\co S^1\times S^3\to BSO(4)$
which is the classifying map for $\pi_S^*\eta$.
In other words, we have the
following commutative diagram
$$
\begin{array}{ccc}
S^1\times S^3&\overset{\tilde{f}}{\longrightarrow}&BSO(4)\\
\downarrow \pi_S&&\downarrow \pi_B\\
S^1\times S^3&\overset{f}{\longrightarrow}&BO(4)\\
\end{array}
$$
By construction, the map $\tilde{f}$ is equivariant under the
action of the group of deck transformations $\mathbb Z/2\mathbb Z$
where the nontrivial element
$i\in\mathbb Z/2\mathbb Z$ acts on  $S^1\times S^3$ by the formula
$(z,q)\mapsto (-z,q)$ and acts on $BSO(4)$ by reversing orientations
of $4$-planes.

Clearly the maps $f$ and $f_0$ are homotopic iff the maps
$\tilde{f}$ and $\tilde{f}_0$ are equivariantly homotopic. By
above we can assume that  $\tilde{f}$ and $\tilde{f}_0$ are
equivariantly homotopic on the 3-skeleton of $S^1\times S^3$. Next
we compute the equivariant cohomology group
$H_{eq}^4(S^1\times S^3,\{\pi_4(BSO(4))\})$ that contains
the obstruction for extending the
homotopy over the $4$-skeleton and show that in our situation the
obstruction has to vanish.

In order to explicitly describe equivariant cochains we have to
identify the action of $\mathbb Z/2\mathbb Z$ on $\pi_4(BSO(4))$.
Recall that $\pi_4(BSO(4))$ classifies the isomorphism classes of orientable
$4$-bundles over $S^4$. On the other hand,
$\pi_4(BSO(4))\cong\pi_3(SO(4))\cong\mathbb Z\oplus\mathbb Z$
where the last isomorphism can be
described explicitly as
$(m,n)\mapsto (q\to q^n\cdot v\cdot q^{m})$
where we identify $\mathbb R^4$
with the quaternions $\mathbb H$ and $S^3$ with the set of unit
quaternions $q$.
According to~\cite{Miln}, the classes
$p_1, e\in H^4(S^4,\mathbb Z)\cong\mathbb Z$
of the bundle $(m,n)$ are given by
$p_1(m,n)=2(m-n)$ and  $e(m,n)=m+n$.
The action of $i$ on $BSO(4)$ sends the canonical oriented
$4$-bundle $\gamma^4$ to $-\g^4$ (i.e the same bundle with
its orientation reversed).
Therefore, $i^*(p_1(\g^4))=p_1(\g^4)$ and $i^*(e(\g^4))=-e(\g^4)$,
and hence the action of $i$ on
$\pi_3(SO(4))$ is given by $i(m,n)=(-n,-m)$.
%

Now once the action is identified, a straightforward computation shows that
$H_{eq}^4(S^1\times S^3,\{\pi_4(BSO(4))\})\cong
(\mathbb Z\oplus\mathbb Z)\slash\mathrm{diagonal}\cong\mathbb Z$.
Let $O_4\in H_{eq}^4(S^1\times S^3,\{\pi_4(BSO(4))\})$
be the obstruction for the equivariant extension of the
homotopy between $\tilde{f}$ and $\tilde{f_0}$ over the
$4$-skeleton.
It remains to show that $O_4$ vanishes.
The double cover $\pi_S$ induces a homomorphism
$$\pi_S^*\co H_{eq}^4(S^1\times S^3,\{\pi_4(BSO(4))\})\to
H_{eq}^4(S^1\times S^3,\pi_S^\#\{\pi_4(BSO(4))\})$$
where the last group is equal to
$H^4(S^1\times S^3,\pi_4(BSO(4))$ because
the pullback bundle of coefficients $\pi_S^\#\{\pi_4(BSO(4))\}$
is trivial.
%
%
We claim that this map is injective.
(Indeed, since both groups are isomorphic to $\mathbb Z$,
it suffices to show that the map is nonzero.
If $\pi_S^*$ were zero, the orientation lift of any nonorientable
$4$-bundle over $S^1\times S^3$ would be trivial
which is certainly not the case since there
exist nonorientable bundles over
$S^1\times S^3$ with nonzero $p_1$. An example
of such a bundle is the Whitney sum of a nontrivial
line bundle and the pullback of a $3$-bundle over $S^4$
with nonzero $p_1$ via a degree one map $S^1\times S^3\to S^4$.)
Since both $f\circ \pi_S$ and $f_0\circ \pi_S$ are null
homotopic we know that $\pi_S^*(O_4)=0$ and
hence $O_4$ vanishes.

\bf{Case of rank three and two.}\rm\
Again, the classifying maps $f$ and $f_0$ are homotopic on the $3$-skeleton
and one has to compute the obstruction
$O_4$ to extending the homotopy over the $4$-skeleton.

If the rank is three, $\mathbb Z/2\mathbb Z$ action on
the coefficient group is trivial and the obstruction group
$H_{eq}^4(S^1\times S^3,\{\pi_4(BSO(3))\})$ reduces
to $H^4(S^1\times S^3,\mathbb Z)$.
As before we have $\pi_S^*(O_4)=0$ and
since in this case $\pi_S^*$ is clearly injective, we
conclude that $O_4$ vanishes.
In the rank two case the obstruction is always zero simply because
$\pi_4(BSO(2))=0$.
\end{proof}

\begin{prop}\label{s1s2}
The total space of any vector bundle over $S^1\times S^2$
has a complete metric of nonnegative curvature such that
the zero section is a soul.
\end{prop}
\begin{proof}
Since all vector bundles over $S^1$ and $S^2$ admit
nonnegatively curved metric such that
the zero sections are souls, it suffices to show that
any vector bundle over $S^1\times S^2$ is isomorphic to the product of
a bundle over $S^1$ and a bundle over $S^2$.

First of all observe that
two rank $k$ vector bundles over $S^1\times S^2$ are isomorphic iff their
restrictions to the two-skeleton $S^1\vee S^2$ are isomorphic.
Indeed, we only need to extend the homotopy
of the classifying maps $S^1\times S^2\to BO(k)$
from $S^1\vee S^2$ to the remaining $3$-cell. This is always
possible since $\pi_3(BO(k))\cong\pi_2(O(n))=0$.

Now let $\xi$ be a vector bundle of rank $k$ over
$S^1\times S^2$ with the classifying map $f\co S^1\times S^2\to BO(k)$.

The case $k=1$ is obvious because line bundles are classified
by $w_1$ and the inclusion $S^1\hookrightarrow S^1\times S^2$
induces an isomorphism on
$H^1(\ , \mathbb Z/2\mathbb Z)$ so that
any line bundle over $S^1\times S^2$ is a pullback
of a bundle over $S^1$.
Similarly, if $k=2$ and $\xi$ is orientable, then
$\xi$ is completely determined by its Euler class.
Since the inclusion $S^2\hookrightarrow S^1\times S^2$
induces an isomorphism on $H^2(\ ,\mathbb Z)$, we conclude that
$\xi$ is a pullback of a bundle over $S^2$.

Assume that $k\ge 3$.
Since $\pi_2(BSO(k))=\mathbb Z/2\mathbb Z$,
there are exactly two $k$-bundles over $S^2$,
namely the trivial bundle and the Whitney sum
of a trivial bundle and a $2$-bundle with nonzero $w_2$.
Let $\kappa$ be a $2$-bundle over $S^2$
that has the same $w_2$ as the restriction of $\xi$
to $S^2$ and let $\lambda$ be a line bundle over $S^1$
with the same $w_1$ as the restriction of $\xi$ to $S^1$.
Finally, let $g$ be the classifying map for the
the Whitney sum of $\kappa\times\lambda$ and
the trivial bundle of rank $(k-3)$.
By construction the restrictions of $f$ and $g$ to
the two-skeleton $S^1\vee S^2$ are homotopic
as needed.

Finally, suppose that $k=2$ and $\xi$ is not orientable.
Note that the orientable two-fold cover $\tilde\xi$
of $\xi$ has zero Euler class.
(Indeed, since $S^2$ represents the generator of
$H_2(S^1\times S^2,\mathbb Z)$
it suffices to show that the intersection number
of $S^2$ and the zero section of $\tilde\xi$ inside the
total space $E(\tilde\xi)$ is zero.
To compute the intersection number put $S^2$ in the general
position to the zero section of $\xi$ and then look at the
the preimage of the manifolds inside $E(\tilde\xi)$.
The covering action of $\mathbb Z/2\mathbb Z$ on $E(\tilde\xi)$
preserves the orientation on the base and changes the orientation
of the total space. Thus, points of intersection come in pairs:
one with plus sign and the other with minus sign. So
the intersection number is zero.)
This implies that the restriction of $\xi$ to $S^2$ has zero
Euler class, and so $\xi|_{S^2}$ is a trivial bundle.
By above $\xi$ is isomorphic
to the product of $\xi|_{S^1}$ and
the rank zero bundle over $S^2$.
\end{proof}

\bibliographystyle{amsalpha}
\bibliography{nn}

\providecommand{\bysame}{\leavevmode\hbox to3em{\hrulefill}\thinspace}
\begin{thebibliography}{GW00}

\bibitem[Bel98]{Bel}
I.~Belegradek, \emph{Pinching, {P}ontrjagin classes, and negatively curved
  vector bundles}, preprint, 1998.

\bibitem[Bro72]{Br}
W.~Browder, \emph{Surgery on simply-connected manifolds}, Springer-Verlag, New
  York, 1972, Ergebnisse der Mathematik und ihrer Grenzgebiete, Band 65.

\bibitem[CG72]{CG}
J.~Cheeger and D.~Gromoll, \emph{On the structure of complete manifolds of
  nonnegative curvature}, Ann. of Math. \textbf{96} (1972), 413--443.

\bibitem[CG72]{CG2}
J.~Cheeger and D.~Gromoll, \emph{The splitting theorem for manifolds of
  nonnegative {R}icci curvature}, J. Differential Geom \textbf{6} (1971/72),
  119--128.

\bibitem[Che73]{Che}
J.~Cheeger, \emph{Some examples of manifolds of nonnegative curvature}, J.
  Differential Geom \textbf{8} (1973), 623--628.

\bibitem[GW00]{GW}
L.~Guijarro and G.~Walschap, \emph{The metric projection onto the soul},
  Trans.~Amer.~Math.~Soc. \textbf{352} (2000), 55--69.

\bibitem[GZ]{GZ2}
K.~Grove and W.~Ziller, \emph{Bundles with nonnegative curvature}, in
  preparation.

\bibitem[GZ99]{GZ}
K.~Grove and W.~Ziller, \emph{Curvature and symmetry of {M}ilnor spheres},
  preprint, 1999.

\bibitem[Hae61]{Hae}
A.~Haefliger, \emph{Plongements diff\'erentiables de vari\'et\'es dans
  vari\'et\'es}, Comment. Math. Helv. \textbf{36} (1961), 47--82.

\bibitem[Hir66]{Hirz}
F.~Hirzebruch, \emph{Topological methods in algebraic geometry},
  Springer-Verlag New York, 1966.

\bibitem[Mil56]{Miln}
J.~Milnor, \emph{On manifolds homeomorphic to the $7$-sphere}, Ann.~of Math.
  \textbf{64} (1956), 399--405.

\bibitem[{\"O}W94]{OW}
M.~{\"O}zaydin and G.~Walschap, \emph{Vector bundles with no soul}, Proc. Amer.
  Math. Soc. \textbf{120} (1994), no.~2, 565--567.

\bibitem[Rig78]{Ri}
A.~Rigas, \emph{{Geodesic spheres as generators of the homotopy groups of ${\rm
  {O}}$, ${\rm {BO}}$}}, J.~Differential Geom. \textbf{13} (1978), no.~4,
  527--545 (1979).

\bibitem[Sie69]{Sieb}
L.~Siebenmann, \emph{On detecting open collars}, Trans.~Amer.~Math.~Soc.
  \textbf{142} (1969), 201--227.

\bibitem[Szc83]{Sz}
A.~Szczepa\'nski, \emph{Aspherical manifolds with the {$\mathbb{Q}$}-homology
  of a sphere}, Mathematika \textbf{30} (1983), 291--294.

\bibitem[Wil98]{Wilk}
B.~Wilking, \emph{On fundamental groups of manifolds of nonnegative curvature},
  to appear in Differential Geom. Appl., 1998.

\bibitem[Yan95]{Yan}
D.~Yang, \emph{On complete metrics of nonnegative curvature on $2$-plane
  bundles}, Pacific J. Math. \textbf{171} (1995), no.~2, 569--583.

\end{thebibliography}

\end{document}